\newcommand{\R}{I\!\!R}
\newcommand{\C}{\mathbb C}

\newcommand{\N}{I\!\!N}
\newcommand{\Id}{\mathbf 1}
\newcommand{\Dim}{\mathrm{dim}}
\newcommand{\Ker}{\mathrm{Ker}}

\newcommand{\Hcal}{\mathcal H}
\newcommand{\Bd}{\mathcal B}
\newcommand{\Bsa}{\mathcal B_{\mathrm{sa}}}

\newcommand{\GLsa}{\mathrm{GL}_{\mathrm{sa}}}
\newcommand{\Spl}{\mathrm{Sp}}
\newcommand{\spl}{\mathrm{sp}}
\newcommand{\dd}{\mathrm d}
\newcommand{\Img}{\mathrm{Im}}

\newcommand{\Conj}{\mathcal K}
\newcommand{\Mconj}{\mathcal K_{\mathrm{m}}}

\newcounter{contador}
\newenvironment{bulletindent}{\setcounter{contador}{0}
\begin{list} {$\bullet$}
{\usecounter{contador}
\setlength{\leftmargin}{10pt}
\setlength{\rightmargin}{10pt}
\setlength{\labelsep}{5pt}
\setlength{\itemsep}{10pt}
\setlength{\topsep}{10pt}}}
{\end{list}}

\documentclass[oneside,draft,letterpaper,11pt]{amsart}

\usepackage{times}
\usepackage{amssymb}

\numberwithin{equation}{section}

\title[Conjugate points on infinite dimensional Riemannian manifolds]{On the Singularities of the Exponential Map in Infinite Dimensional Riemannian Manifolds}
\author[L. Biliotti]{Leonardo Biliotti}
\author[R. Exel]{Ruy Exel}
\author[P.\ Piccione]{Paolo Piccione}
\author[D.\ Tausk]{Daniel V.\ Tausk}
\address{Departamento de Matem\'atica,\hfill\break\indent
Instituto de Matem\'atica e Estat\'\i stica\hfill\break\indent
Universidade de S\~ao Paulo, \hfill\break\indent Caixa Postal 66281,
CEP 05315--970, SP\hfill\break\indent Brazil}
\email{biliotti@math.unifi.it,
exel@mtm.ufsc.br,\hfil\break\indent piccione@ime.usp.br, tausk@ime.usp.br}
\urladdr{http://www.ime.usp.br/\~{}piccione}
\urladdr{http://www.ime.usp.br/\~{}tausk}
\thanks{The last two authors are partially sponsored by CNPq}
\subjclass[2000]{58B20}

\date{November 18th, 2004}

\begin{document}


\swapnumbers

\theoremstyle{plain}\newtheorem*{teo}{Theorem}
\theoremstyle{plain}\newtheorem{prop}{Proposition}
\theoremstyle{plain}\newtheorem{lem}[prop]{Lemma}
\theoremstyle{plain}\newtheorem{cor}[prop]{Corollary}
\theoremstyle{definition}\newtheorem{defin}[prop]{Definition}
\theoremstyle{remark}\newtheorem{rem}[prop]{Remark}
\theoremstyle{plain} \newtheorem{assum}[prop]{Assumption}
\theoremstyle{definition}\newtheorem{example}[prop]{Example}


\begin{abstract}
Using symplectic techniques and spectral analysis of smooth paths of self-adjoint operators, we characterize
the set of conjugate instants along a geodesic in an infinite dimensional Riemannian Hilbert manifold.
\end{abstract}

\maketitle

\begin{section}{Introduction}\label{sec:intro}

Let $(M,\mathfrak g)$ be a (possibly infinite-dimensional) Riemannian Hilbert manifold. For $x\in M$
we denote by $\exp_x$ the exponential map of $(M,\mathfrak g)$ defined in an open subset of $T_xM$.
Let $\gamma:\left[a,b\right[\to M$ ($-\infty<a<b\le+\infty$) be a geodesic, so that $\gamma(t)=\exp_{\gamma(a)}\big((t-a)\gamma'(a)\big)$,
for $t\in\left[a,b\right[$. We denote by $E_t:T_{\gamma(a)}M\to T_{\gamma(t)}M$ the differential
of $\exp_{\gamma(a)}$ at $(t-a)\gamma'(a)$. An instant $t\in\left]a,b\right[$ is said to be {\em conjugate\/} along $\gamma$
if $E_t$ fails to be an isomorphism. Traditionally (see \cite{Grossmann}), a conjugate instant $t$ is called {\em monoconjugate\/}
if $E_t$ fails to be injective and {\em epiconjugate\/} if $E_t$ fails to be surjective; the multiplicity
of a monoconjugate instant $t$ is defined as the dimension of the kernel of $E_t$. As it was already proven in
\cite{Grossmann}, every conjugate instant is epiconjugate (see also Remark~\ref{thm:remEtFredholm}), and it
is convenient to introduce the notion of {\em strictly epiconjugate\/} instant, to denote an instant $t\in\left]a,b\right[$
for which the range of $E_t$ fails to be closed in $T_{\gamma(t)}M$. Unlike in finite-dimensional Riemannian geometry, conjugate instants
along a geodesic can accumulate. The classical example of this phenomenon is given by an infinite dimensional ellipsoid
in $\ell^2$ whose axes form a non discrete subset of the real line (see \cite{Grossmann}); in this example one has
a sequence of monoconjugate instants converging to a strictly epiconjugate instant. The goal of the present article is
to study the distribution of monoconjugate and strictly epiconjugate instants along a geodesic $\gamma$.
We prove the following result:

\begin{teo}
Let $(M,\mathfrak g)$ be a Riemannian Hilbert manifold and let $\gamma:\left[a,b\right[\to M$ ($-\infty<a<b\le+\infty$) be a geodesic; denote
by $\Conj\subset\left]a,b\right[$ the set of conjugate instants and by $\Mconj\subset\Conj$
the set of monoconjugate instants along $\gamma$. Then:
\begin{itemize}
\item[(a)] $\Conj$ is closed in $\left[a,b\right[$;
\item[(b)] the set of strictly epiconjugate instants along $\gamma$ coincides with the set $\Conj'$ of limit points
of $\Conj$, so that $\Conj\setminus\Conj'\subset\Mconj$;
\item[(c)] if $M$ is modeled on a separable Hilbert space then $\Mconj$ is countable.
\end{itemize}
Conversely, given an interval $\left[a,b\right[$, a set $\Conj\subset\left]a,b\right[$ that is closed in
$\left[a,b\right[$, a subset $\Mconj$ of $\Conj$ containing $\Conj\setminus\Conj'$ and a map $\mathfrak m:\Mconj\to\{1,2,\dots,+\infty\}$ then
there exists a conformally flat Riemannian Hilbert manifold $(M,\mathfrak g)$ and a geodesic $\gamma:\left[a,b\right[\to M$
such that the set of conjugate instants along $\gamma$ is $\Conj$, the set of monoconjugate instants along $\gamma$
is $\Mconj$ and the multiplicity of each $t\in\Mconj$ is $\mathfrak m(t)$. Moreover, if $\Mconj$ is countable
then one can choose $M$ to be modeled on a separable Hilbert space.
\end{teo}

The theorem above gives a complete characterization of the
conjugate instants along a geodesic in a Riemannian Hilbert
manifold. Observe that the Theorem implies that if there are no strictly epiconjugate
instants along $\gamma$ (for instance, when the exponential map is
Fredholm) then the set of conjugate instants along any compact
segment of $\gamma$ is finite. This fact had already been proven by Misio\l ek in \cite{Misiolek} under a certain
technical hypothesis that also implies that the exponential map along a geodesic is Fredholm.
We also prove a Morse Index Theorem for geodesics in Riemannian Hilbert manifolds in the case of absence of strictly
epiconjugate instants.

The proof of the several statements in the thesis of our Theorem is scattered along the entire article.
In Section~\ref{sec:symplback} we give a sketchy introduction to the symplectic background necessary to
follow the arguments presented in the article. In Section~\ref{sec:Lagrgeodesic} we show how the study of conjugate instants can be
reduced to the study of curves in the Lagrangian Grassmannian of a Hilbert symplectic space. In Section~\ref{sec:symplsystems}
we give the notion of symplectic system, which are first order differential equations in Hilbert spaces corresponding
to the Jacobi equation along geodesics. In Section~\ref{sec:distribution} we prove statements
(a), (b) and (c) in the thesis of our Theorem. In Section~\ref{sec:xitogamma} we characterize
which curves of Lagrangians arise from Riemannian geodesics and in Section~\ref{sec:spectro} we prove the last
statement in the thesis of our Theorem. Finally, in Section~\ref{sec:Morse}, we prove a Morse Index Theorem.

\end{section}

\begin{section}{Hilbert symplectic spaces and the Lagrangian Grassmannian}
\label{sec:symplback}

A {\em Hilbert symplectic space\/} is a real Hilbert space $(V,\langle\cdot,\cdot\rangle)$ endowed with a
symplectic form; by a {\em symplectic form\/} we mean a skew symmetric bounded bilinear form $\omega:V\times V\to\R$
which is represented by a (anti self-adjoint) isomorphism of $V$. By replacing the inner product of $V$ with an equivalent
one we may always assume that $\omega=\langle J\cdot,\cdot\rangle$, where $J$ is a orthogonal complex structure on $V$.
A subspace $S$ of $V$ is called {\em isotropic\/} if $\omega$ vanishes on $S$, i.e., if $J(S)\subset S^\perp$. A maximal
isotropic subspace of $V$ is called {\em Lagrangian}; a subspace $L\subset V$ is Lagrangian if and only if $J(L)=L^\perp$.
We denote by $\Lambda(V)$ the set of all Lagrangian subspaces of $V$ and we call it the {\em Lagrangian Grassmannian\/}
of $V$.

Let $\Hcal$ be a real Hilbert space; we will denote by $\Bd(\Hcal)$ the space
of all bounded linear operators on $\Hcal$ and by $\Bsa(\Hcal)$ the closed subspace of $\Bd(\Hcal)$
consisting of self-adjoint operators. We consider the orthogonal direct sum $\Hcal^\C=\Hcal\oplus\Hcal$ endowed
with the complex structure $J(x,y)=(-y,x)$.

Given $L_0,L_1\in\Lambda(V)$ with $V=L_0\oplus L_1$,
we consider the map $\varphi_{L_0,L_1}:\mathcal O(L_1)\to\Bsa(L_0)$ defined by
$\varphi_{L_0,L_1}(L)=\rho_{L_1,L_0}T$, where:
\[\mathcal O(L_1)=\big\{L\in\Lambda:V=L\oplus L_1\big\},\]
$T:L_0\to L_1$ is the bounded linear operator whose graph is $L$, and $\rho_{L_1,L_0}:L_1\to L_0$
is the isomorphism given by $\rho_{L_1,L_0}=P_{L_0}J\vert_{L_1}$, where $P_{L_0}$
denotes orthogonal projection onto $L_0$. The maps
$\varphi_{L_0,L_1}$ constitute a smooth atlas on $\Lambda$, so
that $\Lambda$ becomes a smooth Banach manifold. The isomorphism
$\dd\varphi_{L_0,L_1}(L_0):T_{L_0}\Lambda\to\Bsa(L_0)$ is independent
of $L_1\in\mathcal O(L_0)$ and thus for every $L\in\Lambda$ we identify the tangent
space $T_L\Lambda$ with the Banach space $\Bsa(L)$. For $L\in\mathcal O(L_1)$, the differential
of the chart $\dd\varphi_{L_0,L_1}$ at $L$ is given by:
\begin{equation}\label{eq:diffchart}
\dd\varphi_{L_0,L_1}(L)\cdot H=\eta^*H\eta,
\end{equation}
for all $H\in\Bsa(L)$, where $\eta:L_0\to L$ is the isomorphism given by the restriction to $L_0$ of the projection $L\oplus L_1\to L$.

\begin{lem}\label{thm:LvarphiL}
Given Lagrangians $L_0,L_1,L\in\Lambda(V)$ with $L_0,L\in\mathcal O(L_1)$ then:
\begin{itemize}
\item[(a)] $\Ker\big(\varphi_{L_0,L_1}(L)\big)=L\cap L_0$;
\item[(b)] $L+L_0=V$ if and only if $\varphi_{L_0,L_1}(L)$ is surjective;
\item[(c)] $L+L_0$ is closed in $V$ if and only if $\varphi_{L_0,L_1}(L)$ has closed range in $L_0$.
\end{itemize}
\end{lem}
\begin{proof}
Item (a) is trivial. Items (b) and (c) follow from the observation that the isomorphism
$V=L_0\oplus L_1\ni(x,y)\mapsto\big(x,\rho_{L_1,L_0}(y)\big)\in L_0\oplus L_0$ carries
$L+L_0$ to $L_0\oplus\Img\big(\varphi_{L_0,L_1}(L)\big)$.
\end{proof}

\begin{rem}\label{thm:remL0L1Fredholm}
Given Lagrangians $L_0,L_1\in\Lambda(V)$ then $(L_0+L_1)^\perp=J(L_0\cap L_1)$. Thus,
the codimension of the closure of $L_0+L_1$ in $V$ equals the dimension of $L_0\cap L_1$. In particular, $(L_0,L_1)$ is a Fredholm
pair if and only if $\Dim(L_0\cap L_1)<+\infty$ and $L_0+L_1$ is closed in $V$.
\end{rem}

A bounded linear map between between symplectic Hilbert spaces is called a {\em symplectomorphism\/} if it is an isomorphism
that preserves the symplectic forms. The {\em symplectic group\/}
$\Spl(V)$ is the Banach Lie group of symplectomorphisms
of $V$, or equivalently, the group of bounded isomorphisms $T:V\to V$ satisfying the relation $T^*JT=J$, where $T^*$
denotes the adjoint of $T$. The Lie algebra of $\Spl(V)$,
denoted by $\spl(V)$, consists of the bounded linear
operators $X$ on $V$ satisfying $X^*J+JX=0$. An element
$X\in\spl(V)$ can be written in block-matrix form as:
\begin{equation}\label{eq:ABC}
X=\begin{pmatrix} A&B\\C&-A^*
\end{pmatrix},
\end{equation}
where $A\in\Bd(\Hcal)$ and $B,C\in\Bsa(\Hcal)$.

The symplectic group $\Spl(V)$ acts smoothly and transitively on $\Lambda$ and for fixed $L_0\in\Lambda$
the differential of the map $\Spl(V)\ni\Phi\mapsto\Phi(L_0)\in\Lambda$ at a point $\Phi$ in $\Spl(V)$ is given by:
\begin{equation}\label{eq:difbeta}
T_\Phi\Spl(V)=\spl(V)\Phi\ni X\Phi\longmapsto P_{\Phi(L_0)}JX\vert_{\Phi(L_0)}\in\Bsa\big(\Phi(L_0)\big).
\end{equation}

\end{section}

\begin{section}{The curve of Lagrangians associated to a geodesic}
\label{sec:Lagrgeodesic}

Let $(M,\mathfrak g)$ be a
Riemannian Hilbert manifold and let $\gamma:\left[a,b\right[\to M$ be a geodesic. We denote by $R$ the curvature tensor
of $(M,\mathfrak g)$. The
Jacobi equation $\frac{\mathrm D^2v}{\dd t^2}=R(\gamma',v)\gamma'$
along $\gamma$ gives an isomorphism $\Phi_t:(T_{\gamma(a)}M)^\C\to(T_{\gamma(t)}M)^\C$ defined by
$\Phi_t\big(v(a),\frac{\mathrm Dv}{\dd
t}(a)\big)=\big(v(t),\frac{\mathrm Dv}{\dd t}(t)\big)$, for every
Jacobi field $v$ along $\gamma$. The identity:
\[\mathfrak g\big(\tfrac{\mathrm Dv}{\dd t},w\big)-\mathfrak
g\big(v,\tfrac{\mathrm Dw}{\dd t}\big)\equiv\text{constant},\]
satisfied for every Jacobi fields $v$, $w$ along $\gamma$ implies that the map $\Phi_t$ is a symplectomorphism.
The symplectomorphism $\Phi_t$ restricts to a symplectomorphism from
$(\gamma'(a)^\perp)^\C$ to $(\gamma'(t)^\perp)^\C$ and thus it carries Lagrangians of $(\gamma'(a)^\perp)^\C$
to Lagrangians of $(\gamma'(t)^\perp)^\C$. In what follows, we will denote by $V$ the symplectic space
$(\gamma'(a)^\perp)^\C$.
\begin{defin}\label{thm:xigeodesic}
The curve $\xi:\left[a,b\right[\to\Lambda(V)$ defined by $\xi(t)=\Phi_t^{-1}\big(\{0\}\oplus\gamma'(t)^\perp\big)$
is called the {\em curve of Lagrangians associated to the geodesic $\gamma$}.
\end{defin}
We will see later in Corollary~\ref{thm:corxilinha} that $\xi$ is smooth and that $\xi'(t)$ is a negative isomorphism of $\xi(t)$
for all $t\in\left[a,b\right[$.

\begin{rem}\label{thm:remEortog}
For each $t\in\left[a,b\right[$, the differential $E_t$ of the exponential map restricts to a bounded linear
map on $\gamma'(a)^\perp$ taking values in $\gamma'(t)^\perp$; moreover, $E_t\big(\gamma'(a)\big)=\gamma'(t)$.
It follows that the kernel of $E_t$ is contained in $\gamma'(a)^\perp$ and that the range of $E_t$ is closed in $T_{\gamma(t)}M$
if and only if $E_t\big(\gamma'(a)^\perp\big)$ is closed in $\gamma'(t)^\perp$.
We may thus equivalently state the definitions of monoconjugate, epiconjugate, strictly epiconjugate points and
multiplicities in terms of the singularities of the restriction
of $E_t$ to the orthogonal complement of $\gamma$.
\end{rem}

The maps $\Phi_t$ and $E_t$ are related by the identity $(t-a)E_t(x)=P_1\Phi_t(0,x)$, where
$P_1$ denotes the projection onto the first summand of $(\gamma'(t)^\perp)^\C$. This implies that, setting
$L_0=\{0\}\oplus\gamma'(a)^\perp$, we have for all $t\in\left]a,b\right[$:
\begin{itemize}
\item $\xi(t)\cap L_0=\{0\}\oplus\Ker(E_t)$, so that $t$ is monoconjugate along $\gamma$ if and only if $\xi(t)\cap L_0\ne\{0\}$;
\item $\xi(t)+L_0=\Phi_t^{-1}\big(\Img(E_t)\oplus\gamma'(t)^\perp\big)$, so that:
\begin{itemize}
\item[$\diamond$] $t$ is epiconjugate along $\gamma$ if and only if $\xi(t)+L_0\varsubsetneq V$;
\item[$\diamond$] $t$ is strictly epiconjugate along $\gamma$ if and only if $\xi(t)+L_0$ is not closed in $V$.
\end{itemize}
\end{itemize}

\begin{rem}\label{thm:remEtFredholm}
It follows from the items above and from Remark~\ref{thm:remL0L1Fredholm} that the codimension of the closure of the
range of $E_t$ is equal to the dimension of the kernel of $E_t$. In particular, all conjugate instants are epiconjugate
and the operator $E_t$ is Fredholm if and only if $\Dim\big(\Ker(E_t)\big)<+\infty$ and $\Img(E_t)$ is closed.
\end{rem}

\end{section}

\begin{section}{Symplectic systems}
\label{sec:symplsystems}

Let $\Hcal$ be a real Hilbert space. A {\em symplectic system\/} in $\Hcal$ is a smooth curve
$X:\left[a,b\right[\to\spl(\Hcal^\C)$ ($-\infty<a<b\le+\infty$). The curves $A:\left[a,b\right[\to\Bd(\Hcal)$,
$B,C:\left[a,b\right[\to\Bsa(\Hcal)$ corresponding to the block-matrix
decomposition \eqref{eq:ABC} will be called the {\em components\/} of
$X$. The symplectic system is called {\em positive\/} if $B(t)$ is
a positive isomorphism of $\Hcal$ for all
$t\in\left[a,b\right[$. The {\em fundamental solution\/} of the system $X$ is the
smooth curve $\Phi:\left[a,b\right[\to\Spl(\Hcal^\C)$ satisfying
$\Phi'_t=X(t)\Phi_t$ for all $t\in\left[a,b\right[$ and $\Phi_a=\Id$. We
associate a smooth curve $\xi:\left[a,b\right[\to\Lambda(\Hcal^\C)$ to a symplectic
system $X$ by setting $\xi(t)=\Phi_t^{-1}(L_0)$, where
$L_0=\{0\}\oplus\Hcal$. An instant $t\in\left]a,b\right[$ is
called {\em conjugate\/} for $X$ if $\xi(t)\not\in\mathcal
O(L_0)$; we say that $t$ is {\em monoconjugate\/}, {\em epiconjugate}, {\em strictly epiconjugate}
respectively if $L_0\cap\xi(t)\ne\{0\}$, $L_0+\xi(t)\varsubsetneq\Hcal^\C$, $L_0+\xi(t)$ not closed in $\Hcal^\C$.
If $t$ is monoconjugate then the {\em multiplicity\/} of $t$ is the dimension of $L_0\cap\xi(t)$.

\begin{example}\label{exa:geosymplsyst}
Let $(M,\mathfrak g)$ be a Riemannian Hilbert manifold and let $\gamma:\left[a,b\right[\to M$ be a geodesic.
Set $\Hcal=\gamma'(a)^\perp$ and let $R_t:\Hcal\to\Hcal$ be the self-adjoint operator that is conjugated to
the operator $T_{\gamma(t)}M\ni v\mapsto R\big(\gamma'(t),v\big)\gamma'(t)\in T_{\gamma(t)}M$ by parallel transport
along $\gamma$. Then, the Jacobi equation along $\gamma$ can be written as $v''(t)=R_tv(t)$.
We define a symplectic system $X:\left[a,b\right[\to\spl(V)$ in $\Hcal$ by setting
$A(t)=0$, $B(t)=\Id$ and $C(t)=R_t$.
The fundamental solution $\Phi$ of $X$ is given by $\Phi_t\big(v(a),v'(a)\big)=\big(v(t),v'(t)\big)$,
where $v:\left[a,b\right[\to\gamma'(a)^\perp$ is obtained from a Jacobi field along $\gamma$ by parallel transport.
It follows that the curve of Lagrangians $\xi$ associated to the symplectic system $X$ is equal to the
curve of Lagrangians associated to $\gamma$ (see Definition~\ref{thm:xigeodesic}).
\end{example}

A symplectic system $X:\left[a,b\right[\to\Spl(V)$ with components
$A$, $B$, $C$ is called {\em Riemannian\/} if $A(t)=0$ and $B(t)=\Id$ for all $t\in\left[a,b\right[$.

\begin{lem}\label{thm:lemaxilinha}
If $\xi:\left[a,b\right[\to\Lambda(\Hcal^\C)$ is the curve of Lagrangians associated to a symplectic system
$X:\left[a,b\right[\to\spl(\Hcal^\C)$ with components $A$, $B$, $C$ then the derivative of $\xi$ is given by:
\[\langle\xi'(t)x,y\rangle=-\big\langle B(t)\Phi_tx,\Phi_ty\big\rangle,\quad x,y\in\xi(t),\]
for all $t\in\left[a,b\right[$, where $B(t)$ is identified with an element of $\Bsa(L_0)$. In particular, a symplectic
system $X$ is positive if and only if $\xi'(t)$ is a negative isomorphism of $\xi(t)$ for all $t\in\left[a,b\right[$.
\end{lem}
\begin{proof}
It is a simple computation using \eqref{eq:difbeta}, observing that $\xi$ is the composition of the fundamental solution $\Phi$
of $X$ with the map $\Phi\mapsto\Phi(L_0)$.
\end{proof}

\begin{cor}\label{thm:corxilinha}
If $(M,\mathfrak g)$ is a Riemannian Hilbert manifold, $\gamma:\left[a,b\right[\to M$ is a geodesic and
$\xi:\left[a,b\right[\to\Lambda\big((\gamma'(a)^\perp)^\C\big)$ is the curve of Lagrangians associated to $\gamma$
then $\xi$ is smooth and $\xi'(t)$ is a negative isomorphism of $\xi(t)$ for all $t\in\left[a,b\right[$.\qed
\end{cor}

Two symplectic systems $X,\widetilde X:\left[a,b\right[\to\spl(V)$ are
said to be {\em isomorphic\/} if there exists a smooth curve
$\phi:\left[a,b\right[\to\Spl(V)$ with $\phi(t)(L_0)=L_0$ and
$\widetilde\Phi_t=\phi(t)\Phi_t\phi(a)^{-1}$, for all
$t\in\left[a,b\right[$, where $\Phi$ and $\widetilde\Phi$ denote respectively
the fundamental solutions of $X$ and $\widetilde X$. We can write
$\phi(t)$ in block-matrix form as:
\begin{equation}\label{eq:phiZW}
\phi=\begin{pmatrix} Z&0\\{Z^*}^{-1}W&{Z^*}^{-1}
\end{pmatrix},
\end{equation}
where $Z(t)$ is a bounded isomorphism of $\Hcal$ and
$W(t)\in\Bsa(\Hcal)$, for all $t\in\left[a,b\right[$. The components of
isomorphic systems $X$ and $\widetilde X$ are related by the
identities:
\begin{align*}
\widetilde A&=ZAZ^{-1}-ZBWZ^{-1}+Z'Z^{-1},\\
\widetilde B&=ZBZ^*,\\
\widetilde C&={Z^*}^{-1}(WA+C-WBW+A^*W+W')Z^{-1}.
\end{align*}
The curves of Lagrangians $\xi$ and $\tilde\xi$ associated respectively to $X$ and $\widetilde X$ are related by
$\phi(a)\big(\xi(t)\big)=\tilde\xi(t)$, for all $t\in\left[a,b\right[$.

\begin{lem}\label{thm:positRieman}
Every positive symplectic system is isomorphic to a Riemannian
symplectic system.
\end{lem}
\begin{proof}
By setting $W=0$ and $Z=B^{-\frac12}$ in \eqref{eq:phiZW} we see
that every positive symplectic system is isomorphic to one with
$B=\Id$. Let $X$ be a symplectic system with $B=\Id$. We set
$W=\frac12(A+A^*)$ and $Z$ to be the solution of the equation
$Z'=\frac12Z(A^*-A)$ with $Z(a)=\Id$. Then the corresponding
$\phi$ as in \eqref{eq:phiZW} gives the desired isomorphism.
\end{proof}

\end{section}

\begin{section}{The distribution of conjugate instants}
\label{sec:distribution}

In this section we will prove statements (a), (b) and (c) in the thesis of our Theorem. Using the ideas of Sections~\ref{sec:Lagrgeodesic}
and \ref{sec:symplsystems}
and coordinate charts in the Lagrangian Grassmannian we reduce the problem of studying the distribution of conjugate points
along a geodesic to the problem of studying the singularities of a curve of self-adjoint operators on a Hilbert space
having a positive isomorphism as its derivative. This reduction will now be made precise.

Let $(M,\mathfrak g)$ be a Riemannian Hilbert manifold and let $\gamma:\left[a,b\right[\to M$ be a geodesic.
Let $\xi:\left[a,b\right[\to\Lambda(\Hcal^\C)$ denote the curve of Lagrangians associated to $\gamma$, where $\Hcal=\gamma'(a)^\perp$
(see Definition~\ref{thm:xigeodesic}). We have seen in Corollary~\ref{thm:corxilinha} that $\xi$ is smooth and that
$\xi'(t)$ is a negative isomorphism of $\xi(t)$ for all $t$. For any fixed $t_0\in\left[a,b\right[$ we can
find a Lagrangian $L_1\in\mathcal O(L_0)\cap\mathcal O\big(\xi(t_0)\big)$ (see \cite{Compl}) and thus consider a local
representation of $\xi$ around $t_0$ in the chart $\varphi_{L_0,L_1}$; namely, we set:
\begin{equation}\label{eq:Tgeodesic}
T(t)=-\varphi_{L_0,L_1}\big(\xi(t)\big)\in\Bsa(L_0)\cong\Bsa(\Hcal),
\end{equation}
for $t\in\left[a,b\right[$ near $t_0$, so that $T$ is a smooth curve in $\Bsa(\Hcal)$. It follows from \eqref{eq:diffchart}
that $T'(t)$ a positive isomorphism of $\Hcal$ for all $t$. It follows from Lemma~\ref{thm:LvarphiL} that
$t$ is monoconjugate, epiconjugate or strictly epiconjugate along $\gamma$ respectively if $T(t)$ is not injective,
not surjective or $\Img\big(T(t)\big)$ is not closed in $\Hcal$. Moreover, the multiplicity of a monoconjugate instant $t$
is equal to the dimension of $\Ker\big(T(t)\big)$.

In Subsections~\ref{sub:Ruy} and \ref{sub:voltab} we will prove both inclusions involved in statement (b) of our Theorem.
The proof of statement (a) is obtained from the results of Subsection~\ref{sub:voltab} as a simple observation (see Remark~\ref{thm:pontoinicial}).
Finally, in Subsection~\ref{sub:c} we will prove statement (c).

\begin{subsection}{Strictly epiconjugate instants are not isolated}
\label{sub:Ruy}

The fact that strictly epiconjugate instants along a geodesic are not isolated in the set of conjugate instants
is obtained from the following:
\begin{prop}\label{thm:Ruy}
Let $\Hcal$ be a real Hilbert space and let $t\mapsto T(t)\in\Bsa(\Hcal)$ be a differentiable curve around $t_0\in\R$
such that $T'(t_0)$ is a positive isomorphism
of the space $\Hcal$. If the range of $T(t_0)$ is not closed then $t_0$ is a limit point of
the set $\big\{t:\text{$T(t)$ is not invertible}\big\}$.
\end{prop}

The proof of Proposition~\ref{thm:Ruy} requires several preliminary results. We denote by $\sigma(T)$ the
spectrum of a bounded linear operator $T:\Hcal\to\Hcal$.

\begin{lem}\label{thm:LemaSpectro}
Let $T,H\in\Bsa(\Hcal)$ and assume that $\alpha\Id\le H\le\beta\Id$, for some scalars $\alpha,\beta\in\R$. For any
$\lambda\in\sigma(T)$ there exists $\mu\in\sigma(T+H)\cap[\lambda+\alpha,\lambda+\beta]$.
\end{lem}
\begin{proof}
Assume by contradiction that $\sigma(T+H)\cap[\lambda+\alpha,\lambda+\beta]=\emptyset$. Write
$[\lambda+\alpha,\lambda+\beta]=[\lambda_0-r,\lambda_0+r]$. Then $\sigma(T+H-\lambda_0\Id)\cap[-r,r]=\emptyset$, so
$T+H-\lambda_0\Id$ is invertible and $\sigma\big((T+H-\lambda_0\Id)^{-1}\big)\subset\left]\smash{{-\frac1r},\frac1r}\right[$;
thus:
\begin{equation}\label{eq:maiorquer}
\Vert(T+H-\lambda_0\Id)^{-1}\Vert^{-1}>r.
\end{equation}
Since $\alpha\Id\le H\le\beta\Id$, we have $\sigma\big(H+(\lambda-\lambda_0)\Id\big)\subset[-r,r]$.
Hence, by \eqref{eq:maiorquer}:
\[\Vert H+(\lambda-\lambda_0)\Id\Vert\le r<\Vert(T+H-\lambda_0\Id)^{-1}\Vert^{-1};\]
it follows that $T-\lambda\Id=(T+H-\lambda_0\Id)-\big(H+(\lambda-\lambda_0)\Id\big)$
is invertible, contradicting $\lambda\in\sigma(T)$.
\end{proof}

\begin{cor}\label{thm:corKato}
Let $A,B\in\Bsa(\Hcal)$ be fixed. Given $\lambda\in\sigma(A)$, there exists $\mu\in\sigma(B)$
with $\vert\lambda-\mu\vert\le\Vert B-A\Vert$.
\end{cor}
\begin{proof}
Set $T=A$, $H=B-A$, $\alpha=-\Vert B-A\Vert$ and $\beta=\Vert B-A\Vert$ in Lemma~\ref{thm:LemaSpectro}.
\end{proof}

In what follows, we denote by $\GLsa(\Hcal)$ the set of bounded invertible self-adjoint linear operators on $\Hcal$
and by $\mathcal N$ the set:
\[\mathcal N=\big\{A\in\Bsa(\Hcal):\sigma(A)\cap\left]-\infty,0\right[\ne\emptyset\big\}.\]

\begin{cor}\label{thm:CorN}
The following assertions hold:
\begin{itemize}
\item the set $\mathcal N$ is open in $\Bsa(\Hcal)$;
\item the set $\mathcal N\cap\GLsa(\Hcal)$ is closed in $\GLsa(\Hcal)$;
\item the map $\phi:\mathcal N\cap\GLsa(\Hcal)\ni A\mapsto\max\big(\sigma(A)\cap\left]-\infty,0\right[\big)$
is continuous.
\end{itemize}
\end{cor}
\begin{proof}
If $\lambda\in\sigma(A)\cap\left]-\infty,0\right[$ and $B\in\Bsa(\Hcal)$ satisfies $\Vert B-A\Vert<\vert\lambda\vert$ then, by Corollary~\ref{thm:corKato},
$B\in\mathcal N$. This proves that $\mathcal N$ is open in $\Bsa(\Hcal)$.
If $A\in\GLsa(\Hcal)$ is not in $\mathcal N$ then $A\ge c\Id$ for some $c>0$; thus, if $B\in\Bsa(\Hcal)$
satisfies $\Vert B-A\Vert\le\frac c2$ we have $B\ge\frac c2\Id$, which implies $B\not\in\mathcal N$. This proves
that $\mathcal N\cap\GLsa(\Hcal)$ is closed in $\GLsa(\Hcal)$. Let $A\in\mathcal N\cap\GLsa(\Hcal)$ be fixed.
Set $\lambda=\phi(A)$ and let $\varepsilon>0$ be given. Choose $c>0$ with $\sigma(A)\cap[0,c]=\emptyset$.
Given $B\in\mathcal N\cap\GLsa(\Hcal)$ with $\Vert B-A\Vert<\min\{\vert\lambda\vert,\varepsilon,c\}$ we claim
that $\big\vert\phi(B)-\lambda\big\vert<\varepsilon$. Since $\lambda\in\sigma(A)$, Corollary~\ref{thm:corKato}
gives us $\mu\in\sigma(B)$ with $\vert\mu-\lambda\vert<\min\{\vert\lambda\vert,\varepsilon\}$;
thus $\mu<0$ and $\phi(B)\ge\mu>\lambda-\varepsilon$. On the other hand, since $\phi(B)\in\sigma(B)$, Corollary~\ref{thm:corKato}
gives us $\rho\in\sigma(A)$ with $\big\vert\rho-\phi(B)\big\vert<\min\{\varepsilon,c\}$. Thus, since $\sigma(A)\cap[0,c]=\emptyset$,
$\rho$ cannot be positive and hence $\rho\le\lambda$; it follows that $\phi(B)<\rho+\varepsilon\le\lambda+\varepsilon$,
which proves the claim and the continuity of $\phi$.
\end{proof}

\begin{lem}\label{thm:crescelocalafim}
Let $\lambda:I\to\R$ be a continuous map defined on an interval $I\subset\R$ and let $c\in\R$ be fixed.
Assume that for all $x\in I$ there exists $\delta_x>0$ such that:
\begin{equation}\label{eq:lambdac}
\lambda(y)-\lambda(x)\ge c(y-x),
\end{equation}
for all $y\in I\cap\left[x,x+\delta_x\right[$. Then inequality \eqref{eq:lambdac} holds for all $x,y\in I$ with $x\le y$.
\end{lem}
\begin{proof}
The map $I\ni t\mapsto\lambda(t)-ct$ is continuous and locally non decreasing and hence it is
non decreasing on $I$.
\end{proof}

\begin{prop}\label{thm:granderesultado}
Let $T:[a,b]\to\Bsa(\Hcal)$ be a differentiable map with $T'(t)\ge\frac12\Id$ for all $t\in[a,b]$.
If $\lambda_0\in\sigma\big(T(a)\big)\!\cap\left]-\infty,0\right[$ and $b-a\ge2\vert\lambda_0\vert$ then $T(t)$ is not
invertible for some $t\in[a,b]$.
\end{prop}
\begin{proof}
Assume by contradiction that $T(t)\in\GLsa(\Hcal)$ for all $t\in[a,b]$. Since $T(a)\in\mathcal N$,
Corollary~\ref{thm:CorN} implies that $T(t)\in\mathcal N$ for all $t\in[a,b]$ (by connectedness) and that
the map $\lambda:[a,b]\ni t\mapsto\phi\big(T(t)\big)\in\left]-\infty,0\right[$ is continuous. We will show
that, given $\varepsilon>0$, then for all $t\in[a,b]$ there exists $\delta_{t,\varepsilon}>0$ such that:
\begin{equation}\label{eq:lambdast}
\lambda(s)-\lambda(t)\ge\big(\tfrac12-\varepsilon\big)(s-t),
\end{equation}
for all $s\in[a,b]\cap\left[t,t+\delta_{t,\varepsilon}\right[$. If \eqref{eq:lambdast} holds then
Lemma~\ref{thm:crescelocalafim} would give us:
\[\lambda(b)\ge\lambda(a)+\big(\tfrac12-\varepsilon\big)(b-a);\]
taking the limit for $\varepsilon\to0$ we would get:
\[\lambda(b)\ge\lambda(a)+\frac12(b-a)\ge\lambda_0+\vert\lambda_0\vert\ge0,\]
contradicting $\lambda(b)<0$.

To prove \eqref{eq:lambdast}, let $t\in[a,b]$ and $\varepsilon>0$ be fixed. For all $s\in[t,b]$ sufficiently close to $t$,
we have $\Vert T(s)-T(t)-(s-t)T'(t)\Vert\le\varepsilon(s-t)$, which implies:
\[T(s)-T(t)-(s-t)T'(t)\ge-\varepsilon(s-t)\Id.\]
Thus:
\[\Vert T(s)-T(t)\Vert\Id\ge T(s)-T(t)\ge\big(T'(t)-\varepsilon\Id\big)(s-t)\ge\big(\tfrac12-\varepsilon\big)(s-t)\Id,\]
for all $s\in[t,b]$ sufficiently close to $t$. Since $\lambda(t)\in\sigma\big(T(t)\big)$, Lemma~\ref{thm:LemaSpectro}
gives us $\mu\in\sigma\big(T(s)\big)$ with:
\[\lambda(t)+\Vert T(s)-T(t)\Vert\ge\mu\ge\lambda(t)+\big(\tfrac12-\varepsilon\big)(s-t).\]
Thus, for $s\in[t,b]$ sufficiently close to $t$, we have $\mu<0$ and hence:
\[\lambda(s)\ge\mu\ge\lambda(t)+\big(\tfrac12-\varepsilon\big)(s-t),\]
which proves \eqref{eq:lambdast}.
\end{proof}

\begin{lem}\label{thm:somapositivo}
Let $T\in\Bsa(\Hcal)$ and assume that $\sigma(T)\cap\left]c,0\right[=\emptyset$, for some $c<0$. If $H$
is a positive isomorphism of $\Hcal$ with $\Vert H\Vert<\vert c\vert$ then $T+H$ is invertible.
\end{lem}
\begin{proof}
Let $\varepsilon=\min\sigma(H)>0$, so that $-\Vert H\Vert\Id\le-H\le-\varepsilon\Id$.
If $0\in\sigma(T+H)$, Lemma~\ref{thm:LemaSpectro} would give
us $\lambda\in\sigma(T)$ with $c<-\Vert H\Vert\le\lambda\le-\varepsilon<0$,
contradicting our hypothesis on $\sigma(T)$.
\end{proof}

\begin{cor}\label{thm:corPTP}
Let $T\in\Bsa(\Hcal)$ be such that $0$ is a limit point of $\sigma(T)\cap\left]-\infty,0\right[$ and let
$P$ be a positive isomorphism of $\Hcal$. Then $0$ is a limit point of $\sigma(PTP)\cap\left]-\infty,0\right[$.
\end{cor}
\begin{proof}
Assume by contradiction that $\sigma(PTP)\cap\left]c,0\right[=\emptyset$, for some $c<0$. We claim
that $\left]c\Vert P\Vert^{-2},0\right[$ is disjoint from $\sigma(T)$. Given
$\lambda\in\left]c\Vert P\Vert^{-2},0\right[$ then:
\[T-\lambda\Id=P^{-1}(PTP-\lambda P^2)P^{-1},\]
where $-\lambda P^2$ is a positive isomorphism whose norm is less than $\vert c\vert$. The conclusion follows
from Lemma~\ref{thm:somapositivo}.
\end{proof}

\begin{lem}\label{thm:outrogrande}
Let $T:\left[t_0,t_0+\delta\right[\to\Bsa(\Hcal)$ be a differentiable curve with $T'(t_0)$ a positive isomorphism of $\Hcal$.
Assume that $0$ is a limit point of $\sigma\big(T(t_0)\big)\!\cap\left]-\infty,0\right[$. Then $0$ is a limit
point of the set $\big\{t\in\left]t_0,t_0+\delta\right[:\text{$T(t)$ is not invertible}\big\}$.
\end{lem}
\begin{proof}
Replacing $T(t)$ with $T'(t_0)^{-\frac12}T(t)T'(t_0)^{-\frac12}$ and keeping in mind Corollary~\ref{thm:corPTP},
we may assume that $T'(t_0)=\Id$; by continuity we may therefore assume that $T'(t)\ge\frac12\Id$, for all $t$.
Let $\varepsilon\in\left]0,\delta\right[$ be given. We will show that $T(t)$ is not invertible
for some $t\in\left]t_0,t_0+\varepsilon\right]$. Let $\mu\in\sigma\big(T(t_0)\big)\cap\left]{-\frac\varepsilon8},0\right[$
and let $t_1\in\left]0,\frac\varepsilon2\right[$ be such that $\Vert T(t_1)-T(t_0)\Vert<\min\big\{\frac\varepsilon8,\vert\mu\vert\big\}$.
Then, by Corollary~\ref{thm:corKato}, there exists $\lambda_1\in\sigma\big(T(t_1)\big)$ such that
$\vert\mu-\lambda_1\vert<\min\big\{\frac\varepsilon8,\vert\mu\vert\big\}$; thus $-\frac\varepsilon4<\lambda_1<0$.
Note that:
\[\varepsilon-t_1>\frac\varepsilon2>2\vert\lambda_1\vert.\]
Applying Proposition~\ref{thm:granderesultado} to the restriction of $T$ to $[t_1,t_0+\varepsilon]$,
we get that $T(t)$ is not invertible for some $t$ in $[t_1,t_0+\varepsilon]$, which concludes the proof.
\end{proof}

\begin{proof}[Proof of Proposition~\ref{thm:Ruy}]
Since $T(t_0)$ is self-adjoint and its range is not closed, $0$ is a limit point of $\sigma\big(T(t_0)\big)$.
If $0$ is a limit point of $\sigma\big(T(t_0)\big)\!\cap\left]-\infty,0\right[$, the conclusion follows directly from
Lemma~\ref{thm:outrogrande}.
If $0$ is a limit point of $\sigma\big(T(t_0)\big)\!\cap\left]0,+\infty\right[$, apply Lemma~\ref{thm:outrogrande}
to the curve $t\mapsto-T(-t)$.
\end{proof}

\end{subsection}

\begin{subsection}{Non isolated conjugate instants are strictly epiconjugate}
\label{sub:voltab}

The fact that non isolated conjugate instants along a geodesic are strictly epiconjugate is obtained from the following:
\begin{prop}\label{thm:nonisolated}
Let $T:I\to\Bsa(\Hcal)$ be a curve of class $C^1$ defined in some interval $I\subset\R$. Assume that $T(t_0)$
has closed range in $\Hcal$ and that $P_NT'(t_0)\vert_N$
is a positive isomorphism of $N=\Ker\big(T(t_0)\big)$ for some $t_0\in I$ (this is the case, for instance,
if $T'(t_0)$ is a positive isomorphism of $\Hcal$). Then $T(t)$ is an isomorphism of $\Hcal$
for $t\ne t_0$ near $t_0$ in $I$.
\end{prop}

In order to prove Proposition~\ref{thm:nonisolated} we need a preliminary result.

\begin{lem}\label{thm:posneg}
Let $\Hcal$ be a Hilbert space and let $\Hcal=\Hcal_1\oplus\Hcal_2$ be an orthogonal direct
sum decomposition for $\Hcal$ into closed subspaces. Let $A\in\Bsa(\Hcal)$ be such that there exists a constant $c>0$
with $\langle Ax,x\rangle\ge c\langle x,x\rangle$ for all $x\in\Hcal_1$
and $\langle Ay,y\rangle\le-c\langle y,y\rangle$, for all $y\in\Hcal_2$. Then $A$ is an isomorphism of $\Hcal$.
\end{lem}
\begin{proof}
We can write $A$ in block-matrix form as $A=\left(\begin{smallmatrix}A_1&B^*\\B&A_2\end{smallmatrix}\right)$,
where $A_1\in\Bsa(\Hcal_1)$, $A_2\in\Bsa(\Hcal_2)$ and $B:\Hcal_1\to\Hcal_2$ is a bounded linear map. Our hypothesis
tells us that $A_1$ is a positive isomorphism of $\Hcal_1$ and $A_2$ is a negative isomorphism of $\Hcal_2$.
Given $(x,y),(u,v)\in\Hcal_1\oplus\Hcal_2$, the condition $A(x,y)=(u,v)$ is equivalent to $y=A_2^{-1}(v-Bx)$
and:
\[(A_1-B^*A_2^{-1}B)x=u-B^*A_2^{-1}v.\]
Since $A_1$ is a positive isomorphism of $\Hcal_1$ and $-B^*A_2^{-1}B$ is a positive operator on $\Hcal_1$,
then $A_1-B^*A_2^{-1}B$ is a positive isomorphism of $\Hcal_1$. The conclusion follows.
\end{proof}

\begin{proof}[Proof of Proposition~\ref{thm:nonisolated}]
Since $T(t_0)$ is self-adjoint and has closed range, we can find a orthogonal decomposition
$\Hcal=\Hcal_+\oplus\Hcal_-\oplus N$ of $\Hcal$
into closed $T(t_0)$-invariant subspaces
such that $T(t_0)\vert_{\Hcal_+}$ (resp., $T(t_0)\vert_{\Hcal_-}$) is a positive isomorphism of $\Hcal_+$ (resp., a negative
isomorphism of $\Hcal_-$). Obviously there exists a constant $c>0$ such that
$\langle T(t)x,x\rangle\le-c\langle x,x\rangle$ for all $x\in\Hcal_-$,
for $t$ near $t_0$; by Lemma~\ref{thm:posneg}, it suffices to show that for $t\ne t_0$ near $t_0$
there exists $c(t)>0$ with
$\langle T(t)z,z\rangle\ge c(t)\langle z,z\rangle$ for all
$z\in\Hcal_+\oplus N$. Using $T(t)=T(t_0)+T'(t_0)(t-t_0)+\mathrm o(t-t_0)$, we can find
constants $c_1,c_2,c_3>0$ and $\varepsilon>0$ such that, for all $t\in\left]t_0,t_0+\varepsilon\right[\cap I$:
\[\langle T(t)x,x\rangle\ge c_1,\quad\langle T(t)y,y\rangle\ge c_2(t-t_0),\quad
\big\vert\langle T(t)x,y\rangle\big\vert\le c_3(t-t_0),\]
for all $x\in\Hcal_+$, $y\in N$ with $\Vert x\Vert=\Vert y\Vert=1$. Thus, for $t>t_0$ sufficiently near $t_0$:
\[\big\vert\langle T(t)x,y\rangle\big\vert\le\frac12\langle T(t)x,x\rangle^{\frac12}\langle T(t)y,y\rangle^{\frac12},\]
for all $x\in\Hcal_+$, $y\in N$. Hence, for $z=x+y\in\Hcal_+\oplus N$:
\begin{multline*}
\langle T(t)z,z\rangle\ge\langle T(t)x,x\rangle+\langle T(t)y,y\rangle-\langle T(t)x,x\rangle^\frac12
\langle T(t)y,y\rangle^\frac12\\
\ge\frac12\big(\langle T(t)x,x\rangle+\langle T(t)y,y\rangle\big)\ge c(t)\langle z,z\rangle,
\end{multline*}
for $t>t_0$ sufficiently near $t_0$, where $c(t)=\frac12\min\{c_1,c_2(t-t_0)\}$. For $t<t_0$ near $t_0$,
a similar argument shows that $\langle T(t)z,z\rangle\le\bar c(t)\langle z,z\rangle$, for all $z\in\Hcal_-\oplus N$
and some negative constant $\bar c(t)$; we get again from Lemma~\ref{thm:posneg} that $T(t)$ is an isomorphism.
\end{proof}

Using the techniques above we obtain easily a result concerning the variation of the Morse index of a path of self-adjoint operators.
The {\em Morse index\/} of a symmetric bilinear form on a real vector space is the supremum of the dimensions
of all subspaces on which it is negative definite; given an operator $A\in\Bsa(\Hcal)$, we define the {\em Morse index\/}
of $A$ to be the Morse index of the corresponding bilinear form $\langle A\cdot,\cdot\rangle$ on $\Hcal$.

\begin{prop}\label{thm:FredElementary}
Under the assumption of Proposition~\ref{thm:nonisolated}, for $t>t_0$ near $t_0$ in $I$, the Morse index of $T(t)$ is equal
to the Morse index of $T(t_0)$ and for $t<t_0$ near $t_0$ in $I$, the Morse index of $T(t)$ is equal to
the Morse index of $T(t_0)$ plus the dimension of the kernel of $T(t_0)$.
\end{prop}
\begin{proof}
Follows from the proof of Proposition~\ref{thm:nonisolated} and from the following observation:
if $\Hcal=\Hcal_1\oplus\Hcal_2$, with $\langle T(t)\cdot,\cdot\rangle$ positive semi-definite on $\Hcal_1$ and
negative definite on $\Hcal_2$ then the Morse index of $T(t)$ is equal to the dimension of $\Hcal_2$.
\end{proof}

\begin{rem}\label{thm:pontoinicial}
Proposition~\ref{thm:nonisolated} implies in particular that given a geodesic $\gamma:\left[a,b\right[\to M$
then there are no conjugate instants $t$ near $a$. Namely, the curve $T$ in \eqref{eq:Tgeodesic} corresponding
to $\gamma$ near $a$ satisfies $T(a)=0$ and thus Proposition~\ref{thm:nonisolated} implies that $T(t)$
is invertible for $t>a$ near $a$. This observation and the fact that the set of isomorphisms of a Hilbert
space is open in $\Bd(\Hcal)$ implies that the set of conjugate instants along $\gamma$ is closed in
$\left[a,b\right[$.
\end{rem}

\end{subsection}

\begin{subsection}{The set of monoconjugate instants is countable in the separable case}
\label{sub:c}

The fact that the set of monoconjugate instants along a geodesic on a Riemannian Hilbert manifold modeled
on a separable Hilbert space is countable is obtained from the following:
\begin{prop}
Let $\Hcal$ be a separable Hilbert space.
Let $T:I\to\Bsa(\Hcal)$ be a $C^1$ curve with $T'(t)$ a positive isomorphism for all $t\in I$. Then
the set:
\[K_{\mathrm m}=\big\{t\in I:\Ker\big(T(t)\big)\ne\{0\}\big\}\]
is countable.
\end{prop}
\begin{proof}
It suffices to show that the intersection of $K_{\mathrm m}$ with a small neighborhood of some fixed $t_0\in I$
is countable. By replacing $T(t)$ with $T'(t_0)^{-\frac12}T(t)T'(t_0)^{-\frac12}$ we may assume that
$T'(t_0)=\Id$. Applying the mean value inequality to the curve $t\mapsto T(t)-\Id$ we get:
\begin{equation}\label{eq:TtTsmeio}
\Vert T(t)-T(s)-(t-s)\Id\Vert\le\frac12\vert t-s\vert,
\end{equation}
for $t,s\in I$ sufficiently close to $t_0$. For each $t\in K_{\mathrm m}$, choose a unitary vector $v(t)$ in the
kernel of $T(t)$. We have:
\[\vert t-s\vert\;\vert\langle v(t),v(s)\rangle\vert=\big\vert\big\langle\big(T(t)-T(s)-(t-s)\Id\big)v(t),v(s)\big\rangle\big\vert,\]
so that, by inequality \eqref{eq:TtTsmeio}, we have $\big\vert\langle v(t),v(s)\rangle\big\vert\le\frac12$,
for distinct $t,s\in K_{\mathrm m}$ sufficiently close to $t_0$. This implies that $\Vert v(t)-v(s)\Vert\ge1$
and concludes the proof.
\end{proof}

\end{subsection}

\end{section}

\begin{section}{Constructing a geodesic from a curve of Lagrangians}
\label{sec:xitogamma}

Let $\Hcal$ be a real Hilbert space. Let us show now that every Riemannian symplectic system in $\Hcal$ originates from
a Riemannian geodesic:
\begin{lem}\label{thm:Helfer}
Let $R:\left[a,b\right[\to\Bsa(\Hcal)$ be a smooth curve and consider the
Riemannian symplectic system
$X=\left(\begin{smallmatrix}0&\Id\\R&0\end{smallmatrix}\right)$. Consider the
Riemannian manifold $(M,\mathfrak g)$, where $M=\Hcal\oplus\R$,
$\mathfrak g=e^\Omega\mathfrak g_0$, $\mathfrak g_0$ is the inner
product of the orthogonal direct sum $\Hcal\oplus\R$ and
$\Omega:M\to\R$ is defined by $\Omega(x,t)=\langle
R(t)x,x\rangle$, for all\footnote{Here we consider an arbitrary smooth
extension of $R$ to the real line.} $t\in\R$, $x\in\Hcal$. Then
$\gamma:\left[a,b\right[\ni t\mapsto(0,t)\in M$ is a geodesic of
$(M,\mathfrak g)$ and $X$ is the symplectic system corresponding to $(M,\mathfrak g)$
and $\gamma$ as in Example~\ref{exa:geosymplsyst}.
\end{lem}
\begin{proof}
A direct computation shows that the Christoffel symbols of
$(M,\mathfrak g)$ vanish along the axis $\{0\}\times\R$ and thus
$\gamma$ is a geodesic.
The conclusion follows from a straightforward computation
of the curvature tensor of $(M,\mathfrak g)$.
\end{proof}

As in Section~\ref{sec:symplsystems}, we denote by $L_0$ the Lagrangian $\{0\}\oplus\Hcal$.
If $\xi:\left[a,b\right[\to\Lambda$ is a curve with $\xi(a)=L_0$ then a {\em lifting\/} of $\xi$
is a map $\psi:\left[a,b\right[\to\Spl(V)$ with $\psi(a)=\Id$ and $\psi(t)(L_0)=\xi(t)$, for all $t\in\left[a,b\right[$.
\begin{lem}\label{thm:lifting}
Every smooth curve $\xi:\left[a,b\right[\to\Lambda$ with $\xi(a)=L_0$ admits a smooth lifting.
\end{lem}
\begin{proof}
For each $t\in\left[a,b\right[$, set $X(t)=-J\xi'(t)P_{\xi(t)}\in\spl(V)$ and consider the solution
$\Phi$ of the ODE $\Phi'_t=X(t)\Phi_t$ satisfying the initial condition $\Phi_a=\Id$. A simple
computation using \eqref{eq:difbeta} shows that both $t\mapsto\Phi_t(L_0)$ and $\xi$ are integral curves of the time-dependent vector field
$\mathcal V(t)(L)=P_LJX(t)\vert_L\in T_L\Lambda$ on $\Lambda$, both starting at $L_0$. Hence the two curves coincide
and $\Phi$ is a lifting of $\xi$.
\end{proof}

\begin{lem}\label{thm:xiX}
If $\xi:\left[a,b\right[\to\Lambda$ is a smooth curve with $\xi(a)=L_0$ then
there exists a symplectic system $X:\left[a,b\right[\to\spl(V)$ whose associated curve in $\Lambda$
is $\xi$.
\end{lem}
\begin{proof}
By Lemma~\ref{thm:lifting}, we can find a smooth lifting
$\psi:\left[a,b\right[\to\Spl(V)$ of $\xi$. The conclusion is obtained by
setting $X(t)=-\psi(t)^{-1}\psi'(t)$, for all $t$, observing that
the fundamental solution of $X$ is given by
$\Phi_t=\psi(t)^{-1}\psi(a)$.
\end{proof}

\begin{lem}\label{thm:chato1}
Let $(E,\Vert\cdot\Vert)$ be a Banach space, $U\subset E$ be a connected open set, $u\in U$ a fixed point, $\bar\tau:\left[c,b\right[\to U$ a smooth curve and $a\in\R$,
$a<c$. Then there exists $M>0$ such that for all $\eta>0$ there exists a smooth extension
$\tau:\left[a,b\right[\to U$ of $\bar\tau$ with the following properties:
\begin{itemize}
\item $\int_a^c\tau=u(c-a)$;
\item $\big\Vert\tau\vert_{[a,c]}\big\Vert_\infty=\sup_{t\in[a,c]}\Vert\tau(t)\Vert\le M$;
\item $\tau\vert_{[a,c-\eta]}$ is constant.
\end{itemize}
\end{lem}
\begin{proof}
Let $r>0$ be such that the open ball $B(u;r)$ of center $u$ and radius $r$ is contained in $U$ and choose a smooth
curve $\tilde\gamma:\left[c-1,b\right[\to U$ such that $\tilde\gamma(c-1)=u$ and $\tilde\gamma\vert_{\left[c,b\right[}=\bar\tau$. Set
$M=\Vert u\Vert+1+\Vert\tilde\gamma\Vert_\infty$ and choose $\varepsilon>0$ small enough such that $\varepsilon<\eta$
and
\begin{equation}\label{eq:chata}
\frac\varepsilon{c-a-\varepsilon}\Vert\tilde\gamma-u\Vert_\infty<\min\{r,1\}.
\end{equation}
Now, let $\gamma:\left[c-\varepsilon,b\right[\to U$ be a smooth non decreasing reparameterization of $\tilde\gamma$ such that
$\gamma\vert_{\left[c,b\right[}=\bar\tau$ and $\gamma\vert_{[c-\varepsilon,c-\frac\varepsilon2]}\equiv u$. Choose smooth
functions
$\phi_1,\phi_2:\left[a,b\right[\to[0,1]$ with
$\phi_1+\phi_2\equiv1$ and such that the support of $\phi_1$ is contained in $\left[a,c-\frac\varepsilon2\right[$ and
the support of
$\phi_2$ is contained in
$\left]c-\varepsilon,b\right]$. Finally set:
\[\delta=\frac{-\int_a^c\phi_2(\gamma-u)}{\int_a^c\phi_1},\]
and define $\tau=\phi_1(u+\delta)+\phi_2\gamma$. To check that such $\tau$ works observe that $\Vert\delta\Vert$ is
less than or equal to the left hand side of \eqref{eq:chata}.
\end{proof}

\begin{cor}\label{thm:corchato}
Let $\bar\sigma:\left[c,b\right[\to\Bsa(\Hcal)$ be a smooth map such that $\bar\sigma(c)$ and $\bar\sigma'(t)$ are negative isomorphisms
for all $t\in\left[c,b\right[$. Then,
given $a<c$ there exists a smooth extension
$\sigma:\left[a,b\right[\to\Bsa(\Hcal)$ of
$\bar\sigma$ such that
$\sigma(a)=0$, $\sigma(t)$ is a negative isomorphism, for all $t\in\left]a,c\right]$ and $\sigma'(t)$ is a
negative isomorphism for all $t\in\left[a,b\right[$.
\end{cor}
\begin{proof}
Simply apply Lemma~\ref{thm:chato1} to the following objects:
\begin{itemize}
\item $E=\Bsa(\Hcal)$;
\item $U=\big\{B\in\Bsa(\Hcal):\text{$B$ is a negative isomorphism}\big\}$;
\item $u=\dfrac{\bar\sigma(c)}{c-a}$;
\item $\bar\tau=\bar\sigma'$;
\item $\eta>0$ is chosen small enough so that $\eta M<r$, where $r>0$ is such
that the open ball $B\big(\bar\sigma(c);r\big)$ is contained in $U$.
\end{itemize}
Finally,  for $t\in\left[a,b\right[$ define $\sigma(t)=\int_a^t\tau$.
\end{proof}

\begin{prop}\label{thm:extendxi}
Let $\bar\xi:\left[c,b\right[\to\Lambda$ be a smooth curve with $\bar\xi(c)\in\mathcal O(L_0)$ and such that $\bar\xi'(t)$ a negative isomorphism
of $\bar\xi(t)$ for all $t\in\left[c,b\right[$. Then, given $a<c$, there exists a smooth extension $\xi:\left[a,b\right[\to\Lambda$ of $\bar\xi$ with
$\xi(a)=L_0$, $\xi(t)\in\mathcal O(L_0)$, for all $t\in\left]a,c\right[$ and $\xi'(t)$ a negative isomorphism of $\xi(t)$,
for all $t\in\left[a,b\right[$.
\end{prop}
\begin{proof}
Let $L_1\in\Lambda$ be the Lagrangian such that $\varphi_{L_0,\bar\xi(c)}(L_1)$ is the identity of $L_0$; in particular,
$L_1$ is transversal to both $L_0$ and $\bar\xi(c)$. It is easily seen that $\varphi_{L_0,L_1}\big(\bar\xi(c)\big)=-\Id$.
Let $b'\in\left]c,b\right]$ be such that $\bar\xi([c,b'])$ is contained in the domain $\mathcal O(L_1)$ of the chart
$\varphi_{L_0,L_1}$ and define $\bar\sigma:[c,b']\to\Bsa(L_0)\cong\Bsa(\Hcal)$ by
$\bar\sigma=\varphi_{L_0,L_1}\circ\bar\xi\vert_{[c,b']}$. The conclusion follows by an application of
Corollary~\ref{thm:corchato} to $\bar\sigma$, keeping in mind that if $\sigma=\varphi_{L_0,L_1}\circ\xi$ then:
\begin{itemize}
\item $\xi(a)=L_0\Leftrightarrow\sigma(a)=0$;

\item $\xi(t)\in\mathcal O(L_0)\Leftrightarrow\text{$\sigma(t)$ is an isomorphism}$;

\item by formula \eqref{eq:diffchart}, $\xi'(t)$ is a negative isomorphism of $\xi(t)$ if and only if
$\sigma'(t)=\dd\varphi_{L_0,L_1}\big(\xi(t)\big)\cdot\xi'(t)$ is a negative isomorphism of $\Hcal$.\qedhere
\end{itemize}
\end{proof}

\begin{prop}\label{thm:xitogamma}
Let $\xi:\left[c,b\right[\to\Lambda(\Hcal^\C)$ be a smooth curve of Lagrangians such that $\xi'(t)$ is a
negative isomorphism of the space $\xi(t)$ for all $t\in\left[c,b\right[$ and such that $\xi(c)\in\mathcal O(L_0)$.
Then, given $a<c$, there exists a symplectomorphism $\phi\in\Spl(\Hcal^\C)$ with $\phi(L_0)=L_0$ and a Riemannian metric
$\mathfrak g$ on $M=\Hcal\times\R$, conformally equivalent
to the product metric, such that $\gamma:\left[a,b\right[\ni t\mapsto(0,t)\in M$ is a geodesic without conjugate
instants in $\left]a,c\right]$ and such that $\left[c,b\right[\ni t\mapsto\phi\big(\xi(t)\big)\in\Lambda$ is the restriction
to $\left[c,b\right[$ of the curve of Lagrangians associated to $\gamma$.
\end{prop}
\begin{proof}
Extend $\xi$ to $\left[a,b\right[$ as in the statement of Proposition~\ref{thm:extendxi}. Apply Lemma~\ref{thm:xiX}
and get a symplectic system $X$; by Lemma~\ref{thm:lemaxilinha}, $X$ is positive and thus by Lemma~\ref{thm:positRieman},
$X$ is isomorphic to a Riemannian symplectic system $\widetilde X$. The conclusion follows from Lemma~\ref{thm:Helfer}.
\end{proof}

\end{section}

\begin{section}{Constructing a curve of Lagrangians with prescribed singularities}
\label{sec:spectro}

Given a smooth curve $T:\left[c,b\right[\to\Bsa(\Hcal)$ of bounded self-adjoint
operators on a Hilbert space $\Hcal$ with $T(c)$ an isomorphism of $\Hcal$ and $T'(t)$ a positive isomorphism
of $\Hcal$ for all $t\in\left[c,b\right[$, we can produce a smooth curve of Lagrangians $\xi:\left[c,b\right[\to\Lambda(\Hcal^\C)$
satisfying the hypothesis of Proposition~\ref{thm:xitogamma} by setting $\xi(t)=\varphi_{L_0,L_1}^{-1}\big({-T}(t)\big)$, where $L_0=\{0\}\oplus\Hcal$ and
$L_1=\Hcal\oplus\{0\}$ (see also \eqref{eq:diffchart}). Let $(M,\mathfrak g)$ and $\gamma:\left[a,b\right[\to M$
be given by Proposition~\ref{thm:xitogamma}. It follows from Lemma~\ref{thm:LvarphiL} that
$t$ is monoconjugate, epiconjugate or strictly epiconjugate along $\gamma$ respectively if $T(t)$ is not injective,
not surjective or $\Img\big(T(t)\big)$ is not closed in $\Hcal$. Moreover, the multiplicity of a monoconjugate instant $t$
is equal to the dimension of $\Ker\big(T(t)\big)$. In this section we will prove the second part of the statement
of our Theorem. In view of the preceding remarks, it suffices to show the following:
\begin{prop}\label{thm:TKKm}
Let $\Conj\subset\left]c,b\right[$ be closed in $\left[c,b\right[$ and let $\Mconj$ be a subset of $\Conj$ that contains
$\Conj\setminus\Conj'$. Let $\mathfrak m:\Mconj\to\{1,2,\ldots,+\infty\}$ be a map. Then there exists a real Hilbert
space $\Hcal$ and a smooth curve $T:\left[c,b\right[\to\Bsa(\Hcal)$ such that:
\begin{itemize}
\item $T'(t)$ is a positive isomorphism of $\Hcal$ for all $t$;
\item $T(t)$ is not invertible if and only if $t\in\Conj$;
\item $T(t)$ is not injective if and only if $t\in\Mconj$, in which case the dimension of
$\Ker\big(T(t)\big)$ is $\mathfrak m(t)$.
\end{itemize}
If $\Mconj$ is countable we can choose $\Hcal$ to be separable.
\end{prop}

\begin{lem}\label{thm:lemaTKKm}
Let $K$ be a compact subset of the real line and let $K_{\mathrm m}$ be a subset of $K$ that contains
$K\setminus K'$. Let $\mathfrak m:K_{\mathrm m}\to\{1,2,\ldots,+\infty\}$ be a map. Then there exists a real Hilbert
space $\Hcal$ and bounded self-adjoint operator $A\in\Bsa(\Hcal)$ whose spectrum is $K$, whose set of eigenvalues
is $K_{\mathrm m}$ and such that for each $t\in K_{\mathrm m}$, the multiplicity of the eigenvalue $t$ is equal
to $\mathfrak m(t)$. If $K_{\mathrm m}$ is countable we can choose $\Hcal$ to be separable.
\end{lem}
\begin{proof} The proof will be split into steps.
\begin{bulletindent}
\item {\em If $\mathcal X$ is a non empty complete metric space without isolated points then $\mathcal X$ contains
a subset $F$ homeomorphic to the Cantor set $\{0,1\}^{\N}$ and in particular $\mathcal X$ is uncountable and it admits
a Borel probability measure that vanishes on each singleton $\{x\}$}.

One can construct a family of non empty closed sets $F_{\epsilon_1\ldots\epsilon_k}\subset\mathcal X$ indexed on finite sequences
of 0's and 1's such that the diameter of $F_{\epsilon_1\ldots\epsilon_k}$ is less than $\frac1k$,
$F_{\epsilon_1\ldots\epsilon_k}\supset F_{\epsilon_1\ldots\epsilon_k\epsilon_{k+1}}$ and
$F_{\epsilon_1\ldots\epsilon_k0}\cap F_{\epsilon_1\ldots\epsilon_k1}=\emptyset$. Thus, the map $\phi:\{0,1\}^{\N}\to\mathcal X$
such that $\phi(\epsilon_1,\epsilon_2,\ldots)$ is the only point of $\bigcap_{k=1}^\infty F_{\epsilon_1\ldots\epsilon_k}$
is a homeomorphism onto its image $F$. A Borel probability measure on $\mathcal X$ that vanishes on singletons and is
supported on $F$ can be constructed using the classical Cantor--Lebesgue function.

\item {\em If $\mathcal X$ is a complete separable metric space in which all non empty open subsets are uncountable
then there exists a finite Borel measure $\mu$ on $\mathcal X$ that vanishes on singletons and such that
$\mu(U)>0$ for every non empty open subset $U$ of $\mathcal X$}.

Let $(U_n)_{n\ge1}$ be a countable basis of non empty open subsets of $\mathcal X$. Since $U_n$ is a non empty
topologically complete metric space without isolated points, there exists a Borel probability measure $p_n$ on $U_n$
that vanishes on singletons. Set $\mu(A)=\sum_{n=1}^\infty\frac1{2^n}p_n(A\cap U_n)$,
for every Borel subset $A$ of $\mathcal X$.

\item {\em If $\mathcal O$ denotes the union of all countable open subsets of $K$ then $\mathcal O$ is countable and it is
contained in the closure of $K_{\mathrm m}$}.

The set $\mathcal O$ can be written as a countable union of countable open sets and thus $\mathcal O$ is countable. Given $x\in \mathcal O$,
we show that $x$ is in the closure of $K\setminus K'$. For any $\varepsilon>0$ sufficiently small,
$\left]x-\varepsilon,x+\varepsilon\right[\cap K$ is a topologically complete non empty countable metric space
and thus it must contain an isolated point $y$, so that $y\in K\setminus K'$.

\item {\em Conclusion of the proof}.

Let $S$ be a set and $f:S\to\R$ be a map such that $f(S)=K_{\mathrm m}$ and such that for every $x\in K_{\mathrm m}$,
$f^{-1}(x)$ is countable and has cardinality equal to $\mathfrak m(x)$. Set $\Hcal=L^2(K\setminus\mathcal O,\mu)\oplus\ell^2(S)$
and $A=A_0\oplus A_1$, where $A_0\in\Bsa\big(L^2(K\setminus\mathcal O,\mu)\big)$ is the multiplication operator by the
identity of $K\setminus\mathcal O$ and $A_1\in\Bsa\big(\ell^2(S)\big)$ is the multiplication operator by $f$. Note that
$\Hcal$ is separable if and only if $K_{\mathrm m}$ is countable. Thus $\sigma(A)=\sigma(A_0)\cup\sigma(A_1)$.
Clearly, $\sigma(A_1)$ is the closure of $K_{\mathrm m}$. Since every non empty open subset of $K\setminus\mathcal O$ has positive
measure $\mu$ it follows that $\sigma(A_0)=K\setminus\mathcal O$; moreover, since the points of $K\setminus\mathcal O$ have null
measure $\mu$ then $A_0$ has no eigenvalue. The conclusion follows by observing that, since $\mathcal O$ is contained
in the closure of $K_{\mathrm m}$, the compact set $K$ is the union of $\sigma(A_0)$ and $\sigma(A_1)$.\qedhere
\end{bulletindent}
\end{proof}

\begin{proof}[Proof of Proposition~\ref{thm:TKKm}]
If $b<+\infty$, apply Lemma~\ref{thm:lemaTKKm} to $K=\overline\Conj$, $K_{\mathrm m}=\Mconj$ and obtain
a self-adjoint operator $A$ on some Hilbert space $\Hcal$. The desired curve $T$ is obtained by
setting $T(t)=t\cdot\Id-A$, for $t\in\left[c,b\right[$. If $b=+\infty$, choose a strictly increasing smooth
diffeomorphism $\theta:\left[c,b\right[\to\left[0,1\right[$ and apply Lemma~\ref{thm:lemaTKKm} to $K=\overline{\theta(\Conj)}$,
$K_{\mathrm m}=\theta(\Mconj)$, obtaining a self-adjoint operator $A$. The conclusion is obtained by setting
$T(t)=\theta(t)\cdot\Id-A$, for $t\in\left[c,b\right[$.
\end{proof}

\end{section}

\begin{section}{A Morse Index Theorem}
\label{sec:Morse}

In this section we consider a compact segment of geodesic $\gamma:[a,b]\to M$
in a Riemannian Hilbert manifold $(M,\mathfrak g)$. Let $H^1_0(\gamma)$ denote the Hilbert space of vector
fields along $\gamma$ of Sobolev type $H^1$ vanishing at the endpoints endowed with the inner product:
\begin{equation}\label{eq:innprodH1}
\langle v,w\rangle_{H^1}=\int_a^b\mathfrak g\big(\tfrac{\mathrm Dv}{\dd t},\tfrac{\mathrm Dw}{\dd t}\big)\,\dd t.
\end{equation}
We recall that the {\em index form\/} of $\gamma$ is the symmetric bounded bilinear form $I$
on $H^1_0(\gamma)$ defined by:
\[I(v,w)=\int_a^b\mathfrak g\big(\tfrac{\mathrm Dv}{\dd t},\tfrac{\mathrm Dw}{\dd t}\big)
+\mathfrak g\big(R(\gamma'(t),v(t))\gamma'(t),w(t)\big)\,\dd t.\]
For $t\in\left]a,b\right]$, we denote by $I_t:H^1_0\big(\gamma\vert_{[a,t]}\big)\times H^1_0\big(\gamma\vert_{[a,t]}\big)\to\R$
the index form corresponding to the geodesic $\gamma\vert_{[a,t]}$. In what follows we identify the symmetric bilinear
form $I_t$ with the self-adjoint operator on $H^1_0\big(\gamma\vert_{[a,t]}\big)$ that represents it.
We will prove a Morse Index Theorem under the assumption that, for all $t$, the operator $I_t$ has closed range.

While it is very easy to see that the kernel of $I_t$ is isomorphic to the kernel of the differential of the exponential map $E_t$,
it is an interesting point to relate the  range of $I_t$  with the range of $E_t$.
We will first assess this point, aiming at a result relating the closedness of $\Img(I_t)$
and that of $\Img(E_t)$.
In \cite[Theorem 1]{Misiolek} the author gives
sufficient conditions under which $E_t$ is a Fredholm operator. The hypothesis of \cite[Theorem 1]{Misiolek} easily implies
that $I_t$ is a compact perturbation of the identity and hence also a Fredholm operator. We will show
that if $E_t$ has closed range (resp., is Fredholm) then $I_t$ also has closed range (resp., is Fredholm).

In analogy to what was observed in Remark~\ref{thm:remEortog}, it is sufficient to study the restriction of the
index form $I$ to the subspace $H^1_0(\gamma^\perp)$ of $H^1_0(\gamma)$ consisting of vector fields pointwise orthogonal
to $\gamma'$. Namely, $H^1_0(\gamma)$ can be decomposed into the direct sum of $H^1_0(\gamma^\perp)$ and the space
$H^1_0(\gamma^\parallel)$ of vector fields pointwise parallel to $\gamma'$; such decomposition is both $\langle\cdot,\cdot\rangle_{H^1}$-orthogonal
and $I$-orthogonal. Moreover, the restriction of $I$ to $H^1_0(\gamma^\parallel)$ is represented by the identity operator.

Consider the Riemannian symplectic system in $\Hcal=\gamma'(a)^\perp$ associated to the Jacobi
equation along $\gamma$ as described in Example~\ref{exa:geosymplsyst}.
Using parallel transport along $\gamma$ we identify the space $H^1_0(\gamma^\perp)$ with the space $H^1_0\big([a,b],\Hcal\big)$
of maps $v:[a,b]\to\Hcal$ of Sobolev type $H^1$ such that $v(a)=v(b)=0$. The inner product \eqref{eq:innprodH1} is identified with
the inner product $\langle v,w\rangle_{H^1}=\int_a^b\langle v',w'\rangle$ on $H^1_0\big([a,b],\Hcal\big)$ and the
index form $I$ is identified with the bilinear form $I(v,w)=\int_a^b\langle v',w'\rangle+\langle Rv,w\rangle$
on $H^1_0\big([a,b],\Hcal\big)$.

\begin{lem}\label{thm:propimagemI}
The image of $I$ is given by:
\[\Img(I)=\Big\{z\in H^1_0\big([a,b],\Hcal\big):P_1\Phi_b\Big(\int_a^b\Phi_t^{-1}\big(z'(t),0\big)\,\dd t\Big)\in\Img(E_b)\Big\},\]
where $P_1:\Hcal^\C\to\Hcal$ denotes the projection onto the first summand.
\end{lem}
\begin{proof}
Given $v,z\in H^1_0\big([a,b],\Hcal\big)$ then $I(v)=z$ if and only if:
\[\int_a^b\langle v'(t)-z'(t),w'(t)\rangle+\langle R_tv(t),w(t)\rangle\,\dd t=0,\]
for all $w\in H^1_0\big([a,b],\Hcal\big)$. Integration by parts and the fundamental lemma of calculus of variations
lead to the differential equation $(v'-z')'(t)=R_tv(t)$, whose solutions are easily computed as:
\[v(t)=P_1\Phi_t\Big(\int_a^t\Phi_s^{-1}\big(z'(s),0\big)\,\dd s+c\Big),\quad c\in\Hcal^\C.\]
The conclusion follows by taking into account the boundary conditions $v(a)=v(b)=0$.
\end{proof}

\begin{prop}\label{thm:propEtIt}
If $E_t$ has closed range (resp., is Fredholm) then $I_t$ also has closed range (resp., is Fredholm).
\end{prop}
\begin{proof} The operator $H^1_0\big([a,b],\Hcal\big)\ni z\mapsto P_1\Phi_b\Big(\int_a^b\Phi_t^{-1}\big(z'(t),0\big)\,\dd t\Big)
\in\Hcal$ is continuous. The proof of the proposition
follows immediately from Lemma~\ref{thm:propimagemI}, Remark~\ref{thm:remEtFredholm} and the fact that the kernels of $I_t$
and of $E_t$ are isomorphic.
\end{proof}

\begin{rem}
It is plausible that the converse of Proposition~\ref{thm:propEtIt} holds, although the authors are not aware of a proof
of this fact. It is not hard to prove the following alternative characterization for the image of $I$:
\[\Img(I)=\Big\{z\in H^1_0\big([a,b],\Hcal\big):\int_a^bE_t^*R_tz(t)\,\dd t\in\Img(E_b^*)\Big\}.\]
Thus, the converse of Proposition~\ref{thm:propEtIt} holds if the
bounded operator $H^1_0\big([a,b],\Hcal\big)\ni
z\mapsto\int_a^bE_t^*R_tz(t)\,\dd t\in\Hcal$ is surjective. This
is the case, for instance, if the operator
$R\big(\gamma'(t),\cdot\big)\gamma'(t)$ is invertible at some non
conjugate instant $t\in\left]a,b\right]$. We remark also that the
converse of Proposition~\ref{thm:propEtIt} would follow if one
could prove that $I_t$ invertible implies that $t$ is not
conjugate. Namely, if $I_t$ has closed range then the proof of the
Morse Index Theorem below implies that $I_s$ is invertible for
$s\ne t$ near $t$ and thus, if $t$ is conjugate, it should be an
isolated conjugate instant. By our Theorem, this implies that
$E_t$ has closed range.
\end{rem}

\begin{prop}[Morse Index Theorem]
Let $(M,\mathfrak g)$ be a Riemannian Hilbert manifold and $\gamma:[a,b]\to M$ a geodesic.
Assume that for each $t\in\left]a,b\right]$, the self-adjoint operator $I_t$
has closed range. Then, the number of monoconjugate instants along $\gamma$ is finite and
the Morse index of $I=I_b$ is equal to the sum of the (possibly infinite) multiplicities of the monoconjugate
instants $t\in\left]a,b\right[$.
\end{prop}
\begin{proof}
We employ the ideas and some standard computations in \cite[Section~3]{Asian}.
Using affine reparameterization one can identify the space $H^1\big([a,t],\Hcal\big)$ with the space
$H^1\big([0,1],\Hcal\big)$, for all $t\in\left]a,b\right]$. The index form $I_t$ corresponds to a
symmetric bilinear form $\hat I_t$ on $H^1\big([0,1],\Hcal\big)$. We set $\mathcal J_t=(t-a)\hat I_t$,
so that $\mathcal J_a=\lim_{t\to a}\mathcal J_t=\langle\cdot,\cdot\rangle_{H^1}$. Clearly $\mathcal J_t$
has closed range and it depends smoothly on $t\in[a,b]$, moreover, for $t>a$, $\mathcal J_t$ and $I_t$
have the same Morse index. The kernel of $\mathcal J_t$ is isomorphic to the orthogonal complement of the
range of $E_t$ and such isomorphism carries $\frac\dd{\dd t}\mathcal J_t$ to $-\langle\cdot,\cdot\rangle$.
The conclusion follows from Proposition~\ref{thm:FredElementary}
and Remark~\ref{thm:remEtFredholm}.
\end{proof}

As an immediate consequence of the Morse Index Theorem, we prove the infinite dimensional analogue of a classical
result due to Morse and Littauer (see \cite{Warner}).
\begin{cor}
Let $(M,\mathfrak g)$ be a Riemannian Hilbert manifold and assume that $\exp_x$ is a Fredholm map for some $x\in M$.
Then, if $v\in T_xM$ is a singular point of $\exp_x$ (i.e., $\dd\exp_x(v)$ is not an isomorphism) then $\exp_x$
is not injective in any neighborhood of $v$ in $T_xM$.
\end{cor}
\begin{proof}
By the Morse Index Theorem, each conjugate instant produces a jump in the Morse index of the geodesic action
functional and therefore the result is a direct application of Krasnosel'skii's bifurcation theory (see \cite{Krasno}).
The details about the application of bifurcation theory to the study of geodesics can be found in \cite{bifurcacao}.
\end{proof}

\end{section}

\end{document}